\input amstex\documentstyle{amsppt}  
\pagewidth{12.5cm}\pageheight{19cm}\magnification\magstep1
\topmatter
\title From groups to symmetric spaces\endtitle
\author G. Lusztig\endauthor
\address{Department of Mathematics, M.I.T., Cambridge, MA 02139}\endaddress
\thanks{Supported in part by the National Science Foundation}\endthanks
\endtopmatter   
\document

\define\bco{\bar{\co}}

\define\dsv{\dashv}
\define\Lie{\text{\rm Lie }}

\define\frl{\forall}
\define\pe{\perp}
\define\si{\sim}

\define\sqc{\sqcup}

\define\qua{\quad}

\define\hf{\hat f}

\define\bG{\bar G}
\define\bK{\bar K}
\define\bC{\bar C}

\define\lb{\linebreak}

\define\op{\oplus}

\define\part{\partial}
\define\em{\emptyset}

\define\ra{\rangle}
\define\n{\notin}
\define\iy{\infty}
\define\m{\mapsto}
\define\do{\dots}
\define\la{\langle}
\define\bsl{\backslash}

\define\sub{\subset}    
\define\bxt{\boxtimes}
\define\T{\times}
\define\ti{\tilde}
\define\nl{\newline}
\redefine\i{^{-1}}

\define\un{\underline}
\define\ov{\overline}
\define\ot{\otimes}
\define\bbq{\bar{\QQ}_l}

\define\ad{\text{\rm ad}}

\define\Hom{\text{\rm Hom}}
\define\End{\text{\rm End}}

\define\Gal{\text{\rm Gal}}

\define\tr{\text{\rm tr}}

\define\g{\gamma}
\redefine\d{\delta}
\define\e{\epsilon}

\redefine\o{\omega}
\define\p{\pi}
\define\ph{\phi}
\define\ps{\psi}
\define\r{\rho}
\define\s{\sigma}
\redefine\t{\tau}
\define\th{\theta}

\redefine\l{\lambda}
\define\z{\zeta}
\define\x{\xi}

\define\Th{\Theta}

\define\Ph{\Phi}
\define\Ps{\Psi}

\define\kk{\bold k}

\define\CC{\bold C}

\define\FF{\bold F}

\define\NN{\bold N}

\define\QQ{\bold Q}

\define\ZZ{\bold Z}

\define\ca{\Cal A}
\define\cb{\Cal B}

\define\cd{\Cal D}

\define\cf{\Cal F}

\define\ch{\Cal H}

\define\ck{\Cal K}
\define\cl{\Cal L}

\define\cn{\Cal N}
\define\co{\Cal O}
\define\cp{\Cal P}

\define\cs{\Cal S}
\define\ct{\Cal T}

\define\cw{\Cal W}
\define\cz{\Cal Z}

\define\fg{\frak g}

\define\fk{\frak k}

\define\fp{\frak p}

\define\ft{\frak t}

\define\fF{\frak F}

\define\tx{\ti x}

\define\tA{\ti A}
\define\tB{\ti B}

\define\tE{\ti E}

\define\tcp{\ti{\cp}}

\define\BKS{BKS}   
\define\GI{Gi}
\define\KL{KL}
\define\CSI{L1}
\define\FOU{L2}
\define\LM{L3}
\define\AFF{L4}
\define\QG{L5}
\define\ST{S}

This paper is based on a talk given at the conference "Representation theory of
real reductive groups", Salt Lake City, July 2009.

We fix an algebraically closed field $\kk$ of characteristic exponent $p$. (We
assume, except in \S18, that either $p=1$ or $p\gg0$.) We also fix a {\it
symmetric space} that is a triple $(G,\th,K)$ where $G$ is a connected 
reductive algebraic group over $\kk$, $\th:G@>>>G$ is an involution and $K$ is
the identity component of the fixed point set of $\th$ ($K$ is a connected 
reductive algebraic group). We shall often write $(G,K)$ instead of 
$(G,\th,K)$. Let $\fg=\Lie(G),\fk=\Lie(K)$, $\fp=\fg/\fp$. Note that $K$ acts 
naturally on $\fp$ by the adjoint action.

If $H$ is a connected reductive algebraic group over $\kk$ then $H$ gives rise
to a symmetric space $(H\T H,H)$ where $H$ is imbedded in $H\T H$ as the
diagonal; here $\th(a,b)=(b,a)$. (Such a symmetric space is said to be 
{\it diagonal}.)

In this paper we examine various properties/constructions which are known for 
groups (or for diagonal symmetric spaces) and we do some experiments to see to
what extent they generalize to non-diagonal symmetric spaces.

\head Contents \endhead
1. Notation.  

2. Almost diagonal symmetric spaces.

3. A generalization of Schur's lemma.

4. Elliptic curves arising from a symmetric space.

5. Dimension of a nilpotent orbit.

6. Some intersection cohomology sheaves.

7. Generalization of a theorem of Steinberg.

8. $\fF$-thin, $\fF$-thick nilpotent orbits.

9. $\fF$-thin nilpotent orbits and affine canonical bases.

10. Character sheaves on $\fp$.

11. Computation of a Fourier transform.

12. Another computation of a Fourier transform.

13. Computation of a Deligne-Fourier transform.

14. Nilpotent $K$-orbits and conjugacy classes in a Weyl group.

15. Example: $(GL_{2n},Sp_{2n})$.

16. Example: $(SO_{2n},SO_{2n-1})$.

17. Example: $(GL_4,GL_2\T GL_2)$.

18. Final comments.

\subhead 1. Notation\endsubhead 
Let $\cb$ be the variety of Borel subgroups of $G$. 

An element $x\in\fp$ is said to be {\it nilpotent} if the closure of the 
$K$-orbit of $x$ in $\fp$ contains $0$. Let $\cn$ be the set of nilpotent 
elements in $\fp$ (a closed $K$-stable subvariety of $\fp$).

Let $l$ be a prime number invertible in $\kk$ and let $\bbq$ be an algebraic 
closure of $\QQ_l$. Let $\bar{}:\bbq@>>>\bbq$ be the involution obtained by 
transporting complex conjugation under some field isomorphism $\bbq@>>>\CC$.

We denote by $\FF_q$ a finite field with $q$ elements; $\bar{\FF}_q$ denotes an
algebraic closure of $\FF_q$; $\ps:\FF_q@>>>\bbq^*$ denotes a fixed nontrivial 
character.

For a finite set $X$, we denote by $|X|$ the cardinal of $X$.

In the remainder of this subsection we assume that $\kk=\bar{\FF}_q$. For an 
algebraic variety $X$ over $\kk$ let $\cd(X)$ be the bounded derived category 
of constructible $\bbq$-sheaves on $X$. If $f:X@>>>\kk$ is a morphism let 
$\cl_f$ be the inverse image under $f$ of the Artin-Schreier local system on
$\kk$ defined by $f$. 

Let $V$ be an $n$-dimensional vector space over $\kk$ and let $V^*$ be is its 
dual. Let $\la,\ra:V^*\T V@>>>\kk$ be the canonical pairing. Then the local 
system $\cl_{\la,\ra}$ on $V^*\T V$ is well defined. Consider the diagram 
$V^*@<a<<V^*\T V@>b>>V$ where $a$ and $b$ is the first and second projection. 
If $\ck\in\cd(V)$ we set $\fF(\ck)=a_!(b^*(\ck)\ot\cl_{\la,\ra})[n]\in\cd(V^*)$
(Deligne-Fourier transform). If $V$ is endowed with a nondegenerate symmetric 
bilinear form then we can use this to identify $V$ and $V^*$ and to regard 
$\fF(\ck)$ as an object of $\cd(V)$.

\subhead 2. Almost diagonal symmetric spaces\endsubhead
Let $\cz$ be the centre of $G$. Let $(\bG,\bar\th,\bK)$ be the symmetric space
such that $\bG=G/\cz$ and $\bar\th:\bG@>>>\bG$ is the involution induced by 
$\th$. We say that $(\bG,\bar\th,\bK)$ is {\it almost diagonal} (resp. {\it 
quasi-split}) if for any involution $\bar\th':\bG@>>>\bG$ such that $\bar\th$,
$\bar\th'$ induce the same involution of the Weyl group of $\bG$ we have 
$\dim(\bK)\ge\dim(\bK')$ (resp. $\dim(\bK)\le\dim(\bK')$); here $\bK'$ is the 
identity component of the fixed point set of $\bar\th'$. We say that 
$(G,\th,K)$ is almost diagonal (resp. quasi-split) if $(\bG,\bar\th,\bK)$ is 
almost diagonal (resp. quasi-split).

For example, if $(G,K)$ is diagonal then it is almost diagonal and quasi-split.
In addition,

(a) $(GL_{2n},Sp_{2n})$, 

(b) $(SO_{2n},SO_{2n-1}), (n\ge2)$,

(c) $(E_6,F_4)$ 
\nl
are almost diagonal. Also, $(G,G)$ is almost diagonal.

We say that $(G,K)$ is of equal rank if $G,K$ contain a common maximal torus.

\subhead 3. A generalization of Schur's lemma\endsubhead
In this subsection we assume that $\kk=\bar\FF_q$ and that we are given an 
$\FF_q$-rational structure on $G$ compatible with $\th$. Then the finite groups
$G(\FF_q),K(\FF_q)$ are well defined. For any irreducible representation $\r$ 
of $G(\FF_q)$ over $\CC$ let $\r_0$ be the space of $K(\FF_q)$-invariant 
vectors in $\r$.
If $(G,K)=(H\T H,H)$ is diagonal and the $\FF_q$-structure on $G$ comes from an
$\FF_q$-structure on $H$ then $\dim\r_0\in\{0,1\}$. Indeed we have 
$\r=\r'\bxt\r''$ where $\r',\r''$ are irreducible representations of $H(\FF_q)$
and the claim follows from Schur's lemma for $H(\FF_q)$. If $(G,K)$ is as in 
2(a) then again $\dim\r_0\in\{0,1\}$ (see \cite{\BKS}); if $(G,K)$ is as in 
2(b), 2(c) then $\dim\r_0\in\{0,1\}$ (R. Lawther). For general $(G,K)$, it is
not true that $\dim\r_0\in\{0,1\}$; but one can show that $\dim\r_0\le c$ where
$c$ is a constant depending only on $(G,K)$ not on $q$ or $\r$ (see 
\cite{\LM}). 
Hence the algebra of double cosets $\CC[K(\FF_q)\bsl G(\FF_q)/K(\FF_q)]$ is a 
direct sum of matrix algebras of sizes $\le c$. 

\subhead 4. Elliptic curves arising from a symmetric space\endsubhead
Let $\co,\co'$ be two $K$-orbits on $\cb$ and let $gK\in G/K$. Following 
\cite{\GI} we consider the subvariety $\co\cap g\co'$ of $\cb$. This clearly
depends only on the coset $gK$. (In the case where $(G,K)$ is diagonal, the 
variety $\co\cap g\co'$ is a special case of a variety which appears in 
\cite{\CSI, 2.5}.) We will give an example (which I have found around 1990)
where $\co\cap g\co'$ is an elliptic curve with finitely many points removed. 
Let $V$ be a $3$-dimensional $\kk$-vector space with a fixed nondegenerate 
symmetric bilinear form $(,)$. Let $SO(V)$ be the corresponding special 
orthogonal group. We take $(G,K)=(GL(V),SO(V))$. For any subspace $V'$ of $V$ 
we set $V'{}^\pe=\{x\in V;(x,V')=0\}$. Let $Q$ be the set of lines $L$ in $V$ 
such that $L\sub L^\pe$. Let $Q'$ be the set of planes $P$ in $V$ such that 
$P^\pe\sub P$. We identify $\cb$ with the set of pairs $L,P$ where $L$ is a 
line in $V$, $P$ is a plane in $V$ and $L\sub P$. There are four $K$-orbits on
$\cb$:

$\co_0=\{L,P;L\sub P,L\in Q,P\in Q'\}$, $\co_1=\{L,P;L\sub P,L\in Q,P\n Q'\}$,
$\co_2=\{L,P;L\sub P,L\n Q,P\in Q'\}$, $\co_3=\{L,P;L\sub P,L\n Q,P\n Q'\}$.
\nl
The closure of $\co_1$ is $\bco_1=\{L,P;L\sub P,L\in Q\}$. The closure of 
$\co_2$ is $\bco_2=\{L,P;L\sub P,P\in Q'\}$.
We can find $g\in GL(V)$ such that the quadrics $Q,gQ$ in the projective plane
are in general position that is $Q\cap gQ$ consists of four distinct points. 
Let $E=\bco_1\cap g\bco_2=\{L,P;L\sub P,L\in Q,P\in gQ'\}$. The non-singular 
surfaces $\{L,P;L\sub P,L\in Q\}$, $\{L,P;L\sub P,P\in gQ'\}$ in $\cb$ 
intersect transversally hence $E$ is a non-singular curve in $\cb$. For any 
$L\in Q$ let $F_L$ be the fibre of the first projection $E@>>>Q$. Now $F_L$ may
be identified with the set of tangents to the quadric $gQ$ passing through $L$.
It consists of two elements if $L$ is not one of the four points in $Q\cap gQ$
and is a single element otherwise. Thus $E$ is a double covering of the 
projective line $Q$ branched at four points. Hence $E$ is an elliptic curve.
Let $E_{ij}=\bco_i\cap g\bco_j$ ($i,j\in[0,3]$). We have 

$E=E_{00}\sqc E_{02}\sqc E_{10}\sqc E_{12}$. 
\nl
We have 

$E_{00}=\{L,P;L\sub P,L\in Q\cap gQ,P\in Q'\cap gQ'\}=\em$.
\nl
(If $L$ is one of the four elements of $Q\cap gQ$ and $P$ is one of the four 
elements of $Q'\cap gQ'$ then $L\not\sub P$.) We have
$$\align&E_{02}=\{L,P;L\sub P,L\in Q, L\n gQ,P\in Q'\cap gQ'\}\\&=
\{L,P;L\sub P,L\in Q,P\in Q'\cap gQ'\}.\endalign$$
Now for any of the four elements $P\in Q'\cap gQ'$ there is a unique $L\in Q$ 
such that $L\sub P$. Hence $E_{02}$ consists of four points. Similarly, 
$E_{10}$ consists of four points. It follows that
$E_{12}=E-(E_{10}\sqc E_{02})$ is an elliptic curve with eight points removed.

Note that if $(G,K)$ is diagonal then the variety $\co\cap g\co'$ (with
$gK\in G/K$) can never be an elliptic curve with finitely many points removed.

\subhead 5. Dimension of a nilpotent orbit\endsubhead
If $(G,K)$ is diagonal, any $K$-orbit in $\cn$ has even dimension. The same 
holds if $(G,K)$ is almost diagonal, but not in the general case.

Consider for example a $\kk$-vector space $V$ of dimension $N\ge4$ with a fixed
nondegenerate symmetric bilinear form $(,)$ and let $U$ be a codimension one
subspace of $V$ such that $(,)$ is nondegenerate when restricted to $U$. Let 
$SO(V)$, $SO(U)$ be the corresponding special orthogonal groups. We have an 
obvious imbedding $SO(U)\sub SO(V)$ which makes $(SO(V),SO(U))$ into a 
symmetric space. 
In this case $\fp$ may be identified with $U$ with the obvious action of
$SO(U)$. We have $\cn=\{x\in U;(x,x)=0\}$. This is the union of two 
$SO(U)$-orbits $\{0\}$ and $\cn-\{0\}$.
Note that $\dim(\cn-\{0\})=N-2$ is even precisely when $N$ is even that is,
precisely when $(SO(V),SO(U))$ is almost diagonal.

\subhead 6. Some intersection cohomology sheaves\endsubhead
Let $C$ be a $K$-orbit in $\cn$. Let $IC(\bC,\bbq)\in\cd(\fp)$ be the
intersection cohomology complex of the closure $\bC$ of $C$ with coefficients 
in $\bbq$, extended by $0$ on $\fp-\bC$. Let $\ch^iIC(\bC,\bbq)$ be the $i$-th
cohomology sheaf of $IC(\bC,\bbq)$ and let $\ch^i_xIC(\bC,\bbq)$ be its stalk 
at $x\in\fp$.

Assume for example that $(G,K)=(SO(V),SO(U))$, with $V,U,N$ as in \S5. Let 
$C=\cn-\{0\}$. Then $\bC=\cn$. Moreover:

(a) If $N$ is even $\ge4$ then $IC(\bC,\bbq)=\bbq$.

(b) If $N$ is odd $\ge3$ then $\ch^0IC(\bC,\bbq)=\bbq$, 
$\dim\ch_x^{N-3}IC(\bC,\bbq)$ is $1$ if $x=0$ and is $0$ if $x\in\fp-\{0\}$;
$\ch^iIC(\bC,\bbq)=0$ if $i\ne0,N-3$.

Now let $(G,K)=(GL(V'\T V''),GL(V')\T GL(V''))$ where $V',V''$ are $\kk$-vector
spaces of dimension $2$. In this case we may identify $\fp$ with 
$\Hom(V',V'')\T\Hom(V'',V')$ with the obvious action of $K$. Let $C_1\sub\fp$ 
be the set of pairs $A,B$ where $A\in\Hom(V',V''),B\in\Hom(V'',V')$ are such 
that $\ker A=B(V'')$ is $1$-dimensional and $\ker B=A(V')$ is one dimensional.
Note that $C_1$ is a $K$-orbit in $\cn$ of dimension $4$. Moreover $\bC_1$ is 
the set of pairs $A,B$ where $A\in\Hom(V',V''),B\in\Hom(V'',V')$ are such that
$A$ and $B$ are singular and $AB=0,BA=0$. We have

(c) $\ch^0IC(\bC_1,\bbq)=\bbq$, $\dim\ch_x^2IC(\bC_1,\bbq)$ is $2$ if $x=0$ and
is $0$ if $x\in\fp-\{0\}$; $\ch^iIC(\bC_1,\bbq)=0$ if $i\ne0,2$.

\subhead 7. Generalization of a theorem of Steinberg\endsubhead
In this subsection we assume that $\kk=\bar\FF_q$. Assume also that we are 
given an $\FF_q$-rational structure on $G$ compatible with $\th$. Then $\fp$ 
and $\cn$ have a natural $\FF_q$-structure. Steinberg \cite{\ST} has shown that
if $(G,K)$ is diagonal, then the number of elements in $\cn(\FF_q)$ is an even
power of $q$. The same holds if $(G,K)$ is almost diagonal, but not in the 
general case.

For example, if $(G,K)$ is as in 2(a) then $|\cn(\FF_q)|=q^{2n^2-2n}$. (Assume
for simplicity that $G$ is $\FF_q$-split. The $K(\FF_q)$-orbits in $\cn(\FF_q)$
are indexed by the partitions of $n$ in such a way that the stabilizer in 
$K(\FF_q)$ of an element in $\cn(\FF_q)$ has cardinal equal to the cardinal of
the stabilizer in $GL_n(\FF_q)$ of a unipotent element of the corresponding 
type (with $q$ replaced by $q^2$), see \cite{\BKS}. Hence the desired equality
follows from the corresponding equality in the diagonal case.)

If $V,U,(,),N$ are as in \S5 with $V,U,()$ defined over $\FF_q$ and 
$(G,K)=(SO(V),SO(U))$ then $|\cn(\FF_q)|$ is equal to $q^{N-2}$ if $N$ is even
and to $q^{N-2}\pm(q^{(N-1)/2}-q^{(N-3)/2})$ if $N$ is odd. Note that 
$(SO(V),SO(U))$ is almost diagonal precisely when $N$ is even.

\subhead 8. $\fF$-thin, $\fF$-thick nilpotent orbits\endsubhead
In this subsection we assume that $\kk=\bar\FF_q$. Let $C$ be a $K$-orbit in 
$\cn$. Then the Deligne-Fourier transform $\fF(IC(\bC,\bbq))\in\cd(\fp)$ is up
to shift a simple perverse sheaf on $\fp$. (We can choose a 
$K$-invariant non-degenerate symmetric bilinear form on $\fp$.)

We say that $C$ is {\it $\fF$-thin} if the support of $\fF(IC(\bC,\bbq))$ is 
contained in $\cn$ (hence is the closure of a single $K$-orbit in $\cn$, by the
finiteness of the number of $K$-orbits in $\cn$). We 
say that $C$ is {\it $\fF$-thick} if the support of $\fF(IC(\bC,\bbq))$ is 
equal to $\fp$ (hence $\fF(IC(\bC,\bbq))$ restricted to some open dense subset
of $\fp$ is up to shift an irreducible local system).

Note that the $K$-orbit $0$ is always $\fF$-thick. The statement that 

"any $K$-orbit in $\cn$ is $\fF$-thick"
\nl
is true if $(G,K)$ is diagonal. It is also true if $(G,K)$ is as in 2(a) (see
\S14) and if $(G,K)$ is as in 2(b) (in this last case this follows from the 
computations in \S12). It is however false for general $(G,K)$.

For example if $(G,K)=(SO(V),SO(U))$ with $V,U,N$ as in \S5, $N$ odd, and if
$C=\cn-\{0\}$ then $C$ is not $\fF$-thick but $\fF$-thin (this follows from the
computations in \S12). 

If $(G,K)=(GL(V'\T V''),GL(V')\T GL(V''))$ and $C_1$ are as in \S6 then $C_1$ 
is $\fF$-thick (this follows from the computations in \S13) but no $K$-orbit in
$\cn$ other than $0$ and $C_1$ is $\fF$-thick.

From these and other examples it appears that, $\fF$-thin $K$-orbits in $\cn$ 
can exist only if $(G,K)$ is of equal rank.

\subhead 9. $\fF$-thin nilpotent orbits and affine canonical bases\endsubhead
In this subsection we assume that $(G,K)=(GL(V'\T V''),GL(V')\T GL(V''))$ where
$V',V''$ are $\kk$-vector spaces of dimension $N',N''$. In this case we may 
identify $\fp$ with $E:=\Hom(V',V'')\T\Hom(V'',V')$ with the obvious action of
$K$. 

Now $\cn$ consists of all $(A,B)\in E$ such that $AB:V''@>>>V''$ and 
$BA:V'@>>>V'$ are nilpotent. We describe the classification of $K$-orbits in 
$\cn$. This is the same as the (known) classification of isomorphism classes
of nilpotent representations of fixed degree of a cyclic quiver with $2$ 
vertices. Let $\cp=\{1,-1\}\T\ZZ_{>0}$. For $(r,m)\in\cp$ we define 
$g_{r,m}\in\NN\T\NN$ by $g_{r,m}=(m/2,m/2)$ if $m$ is even, 
$g_{r,m}=((m+r)/2,(m-r)/2)$ if $m$ is odd. Let $\tcp=\tcp_{N',N''}$ be the set
of maps $\s:\cp@>>>\NN$ such that $\s$ has finite support and 
$\sum_{(r,m)\in\cp}\s(r,m)g_{r,m}=(N',N'')$. For $\s\in\tcp$ let $\cn_\s$ be 
the set of all $(A,B)\in E$ with the following property: there exists direct 
sum decompositions $V'=\op_{i=1}^sV'_i$, $V''=\op_{i=1}^sV''_i$ such that 

$AV'_i\sub V''_i,BV''_i\sub V'_i$, $\dim(V'_i\op V''_i)=m_i$;

the linear map $T_i:V'_i\op V''_i@>>>V'_i\op V''_i$, $(\x',\x'')\m(B\x'',A\x')$
is nilpotent with a single Jordan block; set $r_i=1$ if $T_i^{m_i-1}V'\ne0$,
$r_i=-1$ if $T_i^{m_i-1}V''\ne0$;

for any $(r,m)\in\cp$, the number of $i\in[1,s]$ such that $(r_i,m_i)=(r,m)$ is
equal to $\s(r,m)$.

Then $(\cn_\s)_{\s\in\tcp}$ are precisely the $K$-orbits on $\cn$.

Let $\tcp^0$ be the set of all $\s\in\tcp$ such that for any $m\ge1$ we have
either $\s(1,m)=0$ or $\s(-1,m)=0$. 

In the remainder of this subsection we assume that $\kk=\bar\FF_q$. We show:

(a) {\it If $\s\in\tcp^0$ then the $K$-orbit $\cn_\s$ is $\fF$-thin.}
\nl
The proof is based on the theory of (affine) canonical bases. From 
\cite{\AFF, 5.9} we see that $\ck=IC(\bar\cn_\s,\bbq)$ is (up to shift) an 
element of the canonical basis attached to the cyclic quiver with $2$ vertices.
Using \cite{\QG, 10.2.3} we see that $\fF\ck$ is (up to shift) also an element
of the canonical basis attached to the same cyclic quiver with the opposite 
orientation. Using again \cite{\AFF, 5.9} we see that 
$\fF\ck=IC(\bar\cn_{\s'},\bbq)$ (up to shift) for a well defined $\s'\in\tcp$ 
(we use an isomorphism of the cyclic quiver with the one with opposed 
orientation). In particular the support of $\fF\ck$ is contained in $\cn$. This
proves (a).

We will show elsewhere that, conversely, if $\s\in\tcp$ and $\cn_\s$ is 
$\fF$-thin, then $\s\in\tcp^0$.

\subhead 10. Character sheaves on $\fp$\endsubhead
In this subsection we assume that $\kk=\bar\FF_q$. 

A simple perverse sheaf on $\fp$ is said to be {\it orbital} if it is 
$K$-equivariant and its support is the closure of a single $K$-orbit in $\cn$.
A simple perverse sheaf on $\fp$ is said to be {\it anti-orbital} (or a {\it 
character sheaf}) if it is of the form $\fF(A)$ for some orbital simple 
perverse sheaf $A$ on $\fp$. (These two definitions are identical to the 
definitions of an orbital or anti-orbital complex on $\fg$ given in 
\cite{\FOU}.) For example, if $C$ is an $\fF$-thin nilpotent $K$-orbit in $\fp$
then $IC(\bC,\bbq)$ is up to shift both an orbital and an antiorbital simple 
perverse sheaf.

\subhead 11. Computation of a Fourier transform\endsubhead
We set $k=\FF_q$. Let $U$ be a $k$-vector space of dimension $M\ge3$ with a 
fixed nondegenerate symmetric bilinear form $(,):U\T U@>>>k$. When $M$ is even
we set $\d=1$ if $(,)$ is split over $k$ and $\d=-1$ otherwise. When $M$ is odd
we set $\d=0$. Define $f:U@>>>\bbq$ by $f(x)=0$ if $(x,x)\ne0$, $f(x)=1$ if 
$(x,x)=0,x\ne0$, $f(x)=1+\d q^{(M-2)/2}$ if $x=0$.

Define $\hf:U@>>>\bbq$ by $\hf(x)=q^{-M/2}\sum_{y\in U}\ps((x,y))f(y)$ for 
$x\in U$ (a Fourier transform). When $M$ is odd we define a function 
$\z:\{x\in U;(x,x)\ne0\}@>>>\{1,-1\}$ by $\z(x)=1$ if $(,)$ is split on te 
subspace $\{x'\in U;(x,x')=0\}$ and $\z(x)=-1$ otherwise. We show:

(a) For $M$ even we have $\hf=f$.

(b) For $M$ odd we have $\hf(x)=\z(x)$ if $(x,x)\ne0$, $\hf(x)=0$ if 
$(x,x)=0,x\ne0$, $\hf(x)=q^{(M-2)/2}$ if $x=0$.
\nl
We have 
$$\align&\hf(0)=q^{-M/2}|\{y\in U;(y,y)=0\}|+\d q^{-1}\\&=
q^{-M/2}(q^{M-1}+\d(q^{M/2}-q^{(M-2)/2}))+\d q^{-1}=q^{(M-2)/2}+\d.\endalign$$
Now assume that $x\in U-\{0\},(x,x)=0$. We can find $x'\in U$ such that
$(x',x')=0$, $(x,x')=1$. Let $U'=\{z\in U;(z,x)=0,(z,x')=0\}$. Any $y\in U$ 
such that $(y,y)=0$ can be written uniquely in the form $y=ax+bx'+z$ where 
$a,b\in k,z\in U'$ and $2ab+(z,z)=0$. We have 
$$\hf(x)=q^{-M/2}\sum_{a,b\in k,z\in U';2ab+(z,z)=0}\ps(b)+\d q^{-1}.\tag c$$
If $b\ne0$ and $z\in U'$ then $a$ in the last sum is determined by 
$a=-(z,z)/(2b)$. Hence
$$\align&\hf(x)=q^{-M/2}\sum_{b\in k^*}|U'|\ps(b)+q^{-M/2}
\sum_{a\in k,z\in U';(z,z)=0}1+\d q^{-1}\\&=-q^{(M-4)/2}+
q^{-(M-2)/2}(q^{M-3}+\d(q^{(M-2)/2}-q^{(M-4)/2}))+\d q^{-1}=\d.\endalign$$
Since $f\m\hf$ is an isometry, we have
$$\align&\sum_{x\in U;(x,x)\ne0}\hf(x)\ov{\hf(x)}=
\sum_{x\in U}f(x)\ov{f(x)}-\sum_{x\in U;(x,x)=0}\hf(x)\ov{\hf(x)}\\&=
|\{x\in U;(x,x)=0,x\ne0\}|+(1+\d q^{(M-2)/2})^2\\&
-|\{x\in U;(x,x)=0,x\ne0\}|\d^2-(q^{(M-2)/2}+\d)^2.\endalign$$
This is $0$ if $M$ is even hence in this case $\hf(x)=0$ for any $x\in U$ such
that $(x,x)\ne0$, proving (a).

In the remainder of this subsection we assume that $M$ is odd. Since $f$ is 
invariant under the $O(U)\T k^*$-action on $V$  ($k^*$ acts by homothety) the 
same holds for $\hf$. Hence there exists $\l,\mu\in\bbq$ such that $\hf(x)=\l$
if $x\in\z\i(1)$, $\hf(x)=\mu$ if $x\in\z\i(-1)$. We fix a $2$-dimensional 
subspace $R$ of $U$ such that $(,)$ is nondegenerate, split on $R$. Let 
$U'=\{x\in U;(x,R)=0\}$. For any $c\in k^*$ we can find $x,x'\in R$ such that 
$(x,x)=c$, $(x',x')=0$, $(x,x')=1$. Any $y\in U$ such that $(y,y)=0$ can be 
written uniquely in the form $y=ax+bx'+z$ where $a,b\in k,z\in U'$ and 
$a^2c+2ab+(z,z)=0$. We have 
$$\hf(x)=\sum_{a,b\in k,z\in U';a^2c+2ab+(z,z)=0}\ps(ac+b).$$
Thus $\hf(x)$ does not depend on $x$ but only on $c$ (and $U'$). We denote it 
by $\ph(c)$. If $a\ne0$ and $z\in U'$ then $b$ in the last sum is determined 
uniquely by $b=-(a^2c+(z,z))/(2a)$. Hence 
$$\align&\ph(c)=\sum_{a\in k^*,z\in U'}\ps(ac-(a^2c+(z,z))/(2a))+
\sum_{b\in k,z\in U';(z,z)=0}\ps(b)\\&
=\sum_{a\in k^*,z\in U'}\ps((ac-(z,z)a\i)/2).\endalign$$
We compute
$$\align&\sum_{c\in k^*}\ph(c)=\sum_{z\in U'}\sum_{a\in k^*}(\ps(-(z,z)a\i)/2)
\sum_{c\in k^*}\ps(ac/2))\\&=-\sum_{z\in U'}\sum_{a\in k^*}\ps(-(z,z)a\i)/2)=
\\&-\sum_{z\in U';(z,z)\ne0}(-1)-\sum_{z\in U';(z,z)=0}(q-1)=(q^M-q^{M-1})
-q^{M-1}(q-1)=0.\endalign$$
Note also that $\ph(c)$ depends only on the image of $c$ in $k^*/k^{*2}$. The 
last sequence of equalities shows that $\ph(c)=-\ph(c')$ if $c,c'\in k^*$,
$c'/c\n k^{*2}$. It follows that $\l=-\mu$.

By a property of Fourier transform we have $\sum_{x\in U}\hf(x)=q^{M/2}f(0)$ 
hence $\sum_{x\in U;(x,x)\ne0}\hf(x)=q^{M/2}-q^{(M-2)/2}$ that is
$$\l|\z\i(1)|+\mu|\z\i(-1)|=|\z\i(1)|-|\z\i(-1)|.$$
Since $\l=-\mu$ and $|\z\i(1)|-|\z\i(-1)|\ne0$ it follows that $\l=1,\mu=-1$.
This proves (b).

\subhead 12. Another computation of a Fourier transform\endsubhead
We set $k=\FF_q$. Let $V',V''$ be two $k$-vector spaces of dimension $2$. Let 
$E=\Hom(V',V'')\T\Hom(V'',V')$. Let $E^\circ$ be the set of all $(A;B)\in E$ 
such that $A$ and $B$ are singular and $AB=0$, $BA=0$. (We write $(A;B)$ 
instead of $(A,B)$ to avoid confusion with an inner product.) Define a 
nondegenerate symmetric bilinear form $(,):E\T E@>>>k$ by
$((A;B),(\tA;\tB))=\tr(A\tB,V'')+\tr(B\tA,V')$. Define $f:E@>>>\bbq$ by 
$f(A;B)=1$ if $(A;B)\in E^\circ-\{0;0\}$, $f(0;0)=1+2q$, $f(A;B)=0$ if 
$(A;B)\in E-E^\circ$. Define $\hf:E@>>>\bbq$ by 
$$\hf(\tA;\tB)=q^{-4}\sum_{(A;B)\in E}\ps(((A;B),(\tA;\tB)))f(A;B).$$
(Fourier transform).

Let $E_{rs}$ be the set of all $(A;B)\in E$ such that $A$ and $B$ are 
nonsingular and $AB:V''@>>>V''$, $BA:V'@>>>V'$ are regular semisimple. Define 
$\z:E_{rs}@>>>\bbq$ by $\z(A:B)=1$ if the two eigenvalues of $AB$ (or $BA$) are
in $F_q^*$; $\z(A:B)=-1$ if the two eigenvalues of $AB$ (or $BA$) are not in 
$F_q^*$. We show that for any $(A;B)\in E_{rs}$ we have
$$\hf(A;B)=q^{-2}\z(A;B).\tag a$$
For $(A:B)\in E$ let $f_1(A:B)$ be the number of $(L';L'')$ where $L'$ is a 
line in $V'$, $L''$ is a line in $V''$ and $BV''\sub L'\sub\ker A$, 
$AV'\sub L''\sub\ker B$. This defines a function $f_1:E@>>>\bbq$. Note that
$f_1(A;B)=f(A;B)$ if $(A;B)\ne(0;0)$ and $f_1(0;0)=(q+1)^2=f(0;0)+q^2$.
Let $\hf_1:E@>>>\bbq$ be the Fourier transform of $f_1$. For $(\tA;\tB)\in E$ 
we have
$$\hf_1(\tA;\tB)=q^{-4}\sum_{L',L''}\sum\Sb (A;B);\\
AV'\sub L'',AL'=0,\\BV''\sub L',BL''=0\endSb\ps(\tr(A\tB,V'')+\tr(B\tA,V'))$$
where $L'$ runs through the lines in $V'$ and $L''$ runs through the lines in 
$V''$. For fixed $L',L''$, the set of $(A;B)$ in the last sum is a $k$-vector
space of dimenion $2$ and $(A,B)\m\tr(A\tB,V'')+\tr(B\tA,V')$ is a linear form
on this vector space which is zero if $\tA L'\sub L''$, $\tB L''\sub L'$ and is
zero otherwise. Thus
$$\hf_1(\tA;\tB)=q^{-2}|\{L',L'';\tA L'\sub L'',\tB L''\sub L'\}|.$$ 
Assume now that $(\tA;\tB)\in E_{rs}$. The condition that $\tA L'\sub L''$,
$\tB L''\sub L'$ is equivalent to $\tB\tA L'=L',\tA L'=L''$. Hence
$$\hf_1(\tA;\tB)=q^{-2}|\{L';\tB\tA L'=L'\}|=q^{-2}(\z(\tA;\tB)+1).$$ 
We have 
$$\hf(\tA;\tB)=\hf_1(\tA;\tB)-q^{-2}=q^{-2}\z(\tA;\tB)$$
and (a) is proved.

\subhead 13. Computation of a Deligne-Fourier transform\endsubhead
The results in this subsection were found around 1990.

Let $V$ be a $\kk$-vector space of dimension $2n$ with a fixed nondegenerate 
symplectic form $\la,\ra:V\T V@>>>\kk$. Let
$$E=\{T\in\End(V);\la T(x),y\ra=\la x,T(y)\ra\qua\frl x,y\in V\}.$$
The non-degenerate symmetric bilinear form $(,):\End(V)\T\End(V)@>>>\kk$ given
by $T,T'\m\tr(TT')$ remains nondegenerate when restricted to $E$. The 
symplectic group $Sp(V)$ acts naturally on $E$, preserving $(,)$. Let $E'$ be 
the set of all $T\in E$ such that $T:V@>>>V$ is semisimple and any eigenspace 
of $T$ is $2$-dimensional. Note that $E'$ is open dense in $E$. Let $E_0$ be 
the set of all $T\in E$ such that $T:V@>>>V$ is nilpotent. Note that $E_0$ is 
$Sp(V)$-stable and the set of $Sp(V)$-orbits on $E_0$ is in natural bijection 
with the set of partitions of $n$. (Any element of $E_0$ has Jordan blocks of 
sizes $n_1,n_1,n_2,n_2,\do,n_t,n_t$ where $n_1\ge n_2\ge\do\ge n_t$ is a 
partition of $n$.)

Let $F$ be the set of all flags $V_*=(V_0\sub V_1\sub V_2\do\sub V_{2n})$ in 
$V$ such that $\dim V_i=i$, $V_{2n-i}=\{x\in V;\la x,V_i\ra=0\}$ for all 
$i\in[0,2n]$. For $T\in E,V_*\in F$ we write $T\dsv V_*$ instead of 
$TV_i\sub V_i$ for all $i$. 

Let $V_*\in F$. Let $E^{V_*}=\{T\in E;T\dsv V_*\}$. Let 
$E_0^{V_*}=\{T\in E_0;T\dsv V_*\}$. We have:

(a) $E^{V_*}$ is exactly the orthogonal of $E_0^{V_*}$ with respect to 
$(,):E\T E@>>>\kk$.
\nl
If $T\in E^{V_*}$, $T'\in E_0^{V_*}$ then $TV_i\sub V_i$, $T'V_i\sub V_{i-1}$ 
for all $i\ge1$. Hence $T'TV_i\sub V_{i-1}$ for all $i$ so that $\tr(T'T)=0$ 
and $(T,T')=0$. Thus $E^{V_*}$ is contained in the orthogonal of $E_0^{V_*}$ 
with respect to $(,):E\T E@>>>\kk$. From the definitions we see that

$\dim E^{V_*}+\dim E_0^{V_*}=n^2+(n^2-n)=2n^2-n=\dim E$
\nl
and (a) follows.

Let $\tE=\{(T,V_*)\in E\T F; T\dsv V_*\}$. Define $\p:\tE@>>>E$ by 
$(T,V_*)\m T$. Let $\ck=\p_!\bbq\in\cd(E)$. Let 
$\tE_0=\{(T,V_*)\in E_0\T F; T\dsv V_*\}$. Define $\p_0:\tE_0@>>>E$ by 
$(T,V_*)\m T$. Let $\ck_0=\p_{0!}\bbq\in\cd(E)$. 

In the remainder of this subsection we assume that $\kk=\bar{\FF}_q$. We show:
$$\fF(\ck_0)\cong\ck[n].\tag b$$
Consider the diagram $E@<a<<E\T E@>b>>E$ where $a$ (resp. $b$) is the first 
(resp. second) projection. By definition we have
$$\fF(\ck_0)=a_!(b^*(\p_{0!}\bbq)\ot\cl_{(,)})[2n^2-n].$$
We have $b^*(\p_{0!}\bbq)=\r_!(\bbq)$ where $\r:E\T\tE_0@>>>E\T E$ is given by
$(T;(T',V_*))\m(T;T')$ and
$$\fF(\ck_0)=a_!(\r_!(\cl_f))[2n^2-n]=\ti\r_!\cl_f[2n^2-n]$$
where $\ti\r=a\r:E\T\tE_0@>>>E$ and $f:E\T\tE_0@>>>\kk$ is given by 
$(T;(T',V_*))\m(T,T')$. We have a partition $E\T\tE_0=Z'\sqc Z''$ where 
$Z'=\{(T;(T',V_*))\in E\T\tE_0;T\dsv V_*\}$,
$Z''=\{(T;(T',V_*))\in E\T\tE_0;T\not\dsv V_*\}$.
Let $\ti\r':Z'@>>>E$, $\ti\r'':Z''@>>>E$ be the restrictions of $\ti\r$; let 
$f':Z'@>>>\kk$, $f'':Z''@>>>E\T E$ be the restrictions of $f$. We have a 
distinguished triangle 
$$(\ti\r''_!\cl_{f''},\ti\r_!\cl_f,\ti\r'_!\cl_{f'})$$
in $\cd(E)$. We show that $\ti\r''_!\cl_{f''}=0$. For each $T\in E$ the fibre 
$Z''_T$ of $\ti\r''$ at $T$ may be identified with 
$\{(T',V_*)\in\tE_0;T\not\dsv V_*\}$ and it is enough to show that
$H^i_c(Z''_T;\cl_{f''})=0$ for all $T\in E$ and all $i$. We can map $Z''_T$ to 
$\{V_*\in F;T\not\dsv V_*\}$ by $(T',V_*)\m V_*$. For each $T\in E$ and
$V_*\in\cf$ such that $T\not\dsv V_*$ let $Z''_{T,V_*}$ be the fibre of the
last map at $V_*$. It is enough to show that
$H^i_c(Z''_{T,V_*};\cl_{f''})=0$ for all $T\in E$, $V_*\in\cf$ such that
$T\not\dsv V_*$ and all $i$. This follows from \cite{\QG, 8.1.13} since
$Z''_{T,V_*}$ is a $\kk$-vector space and the restriction of $f''$ to this
vector space (that is $T'\m(T,T')$) is a nonzero linear form (see (a)).

Thus we have $\ti\r''_!\cl_{f''}=0$. From the distinguished triangle above it
follows that $\ti\r_!\cl_f=\ti\r'_!\cl_{f'}$. From (a) we see that $f'$ is 
identically zero. Hence $\cl_{f'}=\bbq$. Moreover $\ti\r'$ factors as 
$Z'@>s>>\tE@>\p>>E$ where $s(T;(T',V_*))=(T,V_*)$. Thus 
$\ti\r_!\cl_f=\p_!s_!\bbq$. Note that $s$ is a vector bundle; its fibre over 
$(T,V_*)$ may be identified with $\{T'\in E_0;T'\dsv V_*\}$, a vector space of 
dimension $n^2-n$. Hence $s_!\bbq=\bbq[-2n^2+2n]$ (we ignore the Tate twist).
This proves (b).

Let $\ca_0$ be the set of isomorphism classes of simple perverse sheaves on $E$
which appear (possibly with a shift) as a direct summand of $\ck_0$. Let $\ca$
be the set of isomorphism classes of simple perverse sheaves on $E$ which 
appear (possibly with a shift) as a direct summand of $\ck$. From (b) we see 
that

(c) $\fF$ defines a bijection $\ca_0@>\si>>\ca$.
\nl
Let $\ca'_0$ be the set of isomorphism classes of $Sp(V)$-equivariant simple 
perverse sheaves on $E$ with support contained in $E_0$. Clearly 

(d) $\ca_0\sub\ca'_0$.
\nl
Let $\p':\p\i(E')@>>>E'$ be the restriction of $\p$. Let $\ca'$ be the set of 
isomorphism classes of simple perverse sheaves on $E$ which appear (possibly 
with a shift) as a direct summand of $\p'_!\bbq$. Clearly we have a natural 
injective map

(e) $\ca'@>>>\ca$.
\nl
From (c),(d),(e) we see that 

(f) $|\ca'|\le|\ca|=|\ca_0|\le|\ca'_0|$.
\nl
Now $\p'$ is a composition $\p\i(E')@>\t>>E''@>\s>>E'$ where $E''$ is the set 
of pairs $(T,\o)$ where $T\in E'$ and $\o$ is an indexing 
$E_{\o(1)},E_{\o(2)},\do,E_{\o(n)}$ 
of the eigenspaces of $T$ by $[1,n]$; $\s$ is defined by $(T,\o)\m T$;
$\t$ is given by $(T,V_*)\m(T,\o)$ with $\o$ defined as follows: for
$i\in[1,n]$ we have $V_i=L_1+L_2+\do+L_i$ where $L_1,L_2,\do L_i$ are lines in
$E_{\o(1)},E_{\o(2)},\do,E_{\o(i)}$ respectively.
Note that $\t$ is a fibre bundle with fibres isomorphic to a product of $n$ 
projective lines. Moreover $\s$ is a finite principal covering with group 
$S_n$, the symmetric group in $n$ letters. It follows that $\t(\bbq)$ is a 
direct sum of shifts of $\bbq$. Hence the objects of $\ca'$ are (up to shift) 
the same as the direct summands of $\s_!\bbq$. Now $E''$ is irreducible (since
$E'$ is so). It follows that the direct summands of $\s_!\bbq$ (up to shift and
isomorphism) are in natural bijection with the irreducible representations of 
$S_n$ (up to isomorphism). Thus, $|\ca'|=p_n$ the number of partitions of $n$.
Now the objects of $\ca'_0$ are in bijection with the $Sp(V)$-orbits on $E_0$
(the stabilizer in $Sp(V)$ of an element in $E_0$ is connected). Hence 
$|\ca'_0|=p_n$. Using now (f) we see that $|\ca'|=|\ca|=|\ca_0|=|\ca'_0|=p_n$.
In particular we see that any object in $\ca'_0$ is in $\ca_0$, any object in 
$\ca$ has support equal to $E$ and $\fF$ defines a bijection between $\ca'_0$ 
and $\ca$. Thus, 

(g) if $A$ is an $Sp(V)$-equivariant simple perverse sheaf on $E$ with support 
contained in $E_0$, then $\fF(A)$ has support equal to $E$.

\subhead 14. Nilpotent $K$-orbits and conjugacy classes in a Weyl group
\endsubhead
In this subsection we assume that $p=1$. We show how the construction of 
\cite{\KL, 9.1} extends to the case of symmetric spaces.

We identify $\fp$ with $\{x\in\fg;\th(x)=-x\}$ in the obvious way. Let $\ct$ 
be the set of subspaces $\ft$ of $\fp$ such that for some torus $T$ in $G$ the
Lie algebra of $T$ equals $\ft$ and such that $\ft$ has maximum possible
dimension. It is known that $K$ acts transitively on $\ct$ by conjugation. For
$\ft\in\ct$ let $\cw$ be the group of components of the normalizer of $\ft$ in
$K$. Let $\un{\cw}$ be the set of conjugay classes in the finite group $\cw$; 
this is independent of the choice of $\ft$.

Let $\Ph=\kk((\e))$, $\ph=\kk[[\e]]$ where $\e$ is an indeterminate. Let 
$\fg_\Ph=\Ph\ot\fp$, $\fp_\Ph=\Ph\ot\fp$, $\fp_\ph=\ph\ot\fp$. Also, the 
groups $G(\Ph),K(\Ph)$ are well defined and the set $\ct(\Ph)$ of $\Ph$-points
of $\ct$ is well defined (it is a set of subspaces of $\fp_\Ph$). 

The group $K(\Ph)$ acts naturally by conjugation on $\ct(\Ph)$; as in 
\cite{\KL, \S1, Lemma 2} we see that the set of $K(\Ph)$-orbits on $\ct(\Ph)$ 
is naturally in $1-1$ correspondence with the set $\un{\cw}$. For 
$\g\in\un{\cw}$ let $\co_\g$ be the $K(\Ph)$-orbit on $\ct(\Ph)$ corresponding
to $\g$.

Let $\fp_{\Ph,rs}$ be the set of all $\x\in\fp_\Ph$ for which there is a unique
$\ft'\in\ct(\Ph)$ such that $\x\in\ft'$; we then set $\ft'_\x=\ft'$. For 
$\g\in\un{\cw}$ let $\fp_{\Ph,rs,\g}$ be the set of all $\x\in\fp_{\Ph,rs}$ 
such that $\ft'_\x\in\co_\g$.

Let $x\in\cn$. As in \cite{\KL, 9.1} there is a unique element $\g\in\un{\cw}$
such that the following holds: there exists a "Zariski open dense" subset $V$ 
of $x+\e\fp_\ph$ (a subset of $\fp_\ph$) such that $V\sub\fp_{\Ph,rs,\g}$. We 
set $\Ps(x)=\g$. Now $\Ps$ is a map $\cn@>>>\un{\cw}$. This map is constant on
$K$-orbits hence it defines a map from the set of $K$-orbits on $\cn$ to 
$\un{\cw}$ denoted again by $\Ps$. (In the case where $(G,K)=(H\T H,H)$ is 
diagonal, $\Ps$ reduces to the map defined in \cite{\KL, 9.1}.)

If $(G,K)$ is diagonal then $\Ps$ is expected to be injective (this is known to
be true in almost all cases). If $(G,K)$ is as in 2(a) then $\Ps$ is bijective
(see \S15). If $(G,K)$ is as in 2(b) then $\Ps$ is bijective (see \S16). For 
general $(G,K)$, $\Ps$ is neither injective nor surjective (see \S17).

\subhead 15. Example: $(GL_{2n},Sp_{2n})$\endsubhead 
In this subsection we assume that 
$(G,K)=(GL(V),Sp(V))$ with $\kk=\CC$ and we keep the notation of \S13. 
Let $\Ph,\ph$ be as in \S14. Let $\bar\Ph$ be an algebraic closure of $\Ph$. 
Define $v:\bar\Ph@>>>\QQ\cup\{\iy\}$ by $v(0)=\iy$ and 
$v(a_0v^f+\text{higher powers of }v)=f$ if $a_0\in\CC^*$, $f\in\QQ$.
Let $E_\Ph=\Ph\ot E,E_\ph=\ph\ot E$, $V_{\bar\Ph}=\bar\Ph\ot V$. Let $E'_\Ph$ 
be the set of all $T\in E_\Ph$ such that $T$ defines a semisimple endomorphism
of $V_{\bar\Ph}$ whose eigenspaces are all $2$-dimensional. For any conjugacy 
class $\s$ in the symmetric group in $n$ letters let $E'_{\Ph,\s}$ be the set 
of all $T\in E'_\Ph$ such that, if the eigenvalues of $T$ are denoted by 
$\l_1=\l_2,\l_3=\l_4,\do,\l_{2n-1}=\l_{2n}$ (elements of $\CC((\e^{1/N}))$ for
some $N\ge1$) then the element $\g$ of $\Gal(\CC((\e^{1/N}))/\Ph)$ which maps 
$\e^{1/N}$ to $\exp(2\pi\sqrt{-1}/N)\e^{1/N}$ permutes $\l_2,\l_4,\do,\l_{2n}$
according to a permutation in $\s$. Note that $E'_\Ph=\sqc_\s E'_{\Ph,\s}$.

Let $x\in E_0$. Recall that $x:V@>>>V$ has Jordan blocks of sizes 
$$n_1,n_1,n_2,n_2,\do,n_t,n_t$$ 
where $n_1\ge n_2\ge\do\ge n_t$ is a partition of $n$. The following statement
is easily verified (it can be reduced to the case where $t=1$). 

(a) There exists $J\in E$ such that $x+\e J:V_\Ph@>>>V_\Ph$ has eigenvalues

$\exp(2\p\sqrt{-1}s/n_i)a_i\e^{1/n_i}\in\bar\Ph$
\nl
$(i=1,2,\do,2t, s=1,2,\do,n_i)$ where $a_2,a_4,\do,a_{2t}$ are distinct in 
$\CC$ and $a_1=a_2,a_3=a_4,\do,a_{2t-1}=a_{2t}$. In particular we have
$x+\e J\in E'_\Ph$.

Now let $X=X_0+\e X_1+\e^2 X_2+\do\in E_\ph$ where $X_0,X_1,\do$ are in $E$.
Let $\tx=x+\e X\in E_\ph$. Let $\l_1,\l_2,\do,\l_{2n}$ be the eigenvalues of 
$\tx$ in $\bar\Ph$. We can assume that $\l_1=\l_2,\l_3=\l_4,\do,
\l_{2n-1}=\l_{2n}$ and $v(\l_2)\le v(\l_4)\le\do\le v(\l_{2n})$.
Let $\mu_1,\mu_2,\do,\mu_{2n}$ be the eigenvalues of $x+\e J$. We can assume 
that $\mu_1=\mu_2,\mu_3=\mu_4,\do,\mu_{2n-1}=\mu_{2n}$ and 
$v(\mu_2)\le v(\mu_4)\le\do\le v(\mu_{2n})$. Applying \cite{\KL, 9.4} we have 
$v(\l_1\l_2\do\l_s)\ge v(\mu_1\mu_2\do\mu_s)$ for $s=1,2,\do,2n$. Using this 
for $s=2,4,\do,2n-2$ and adding the resulting inequalities we deduce
$$v(\l_2^{2n-2}\l_4^{2n-4}\do\l_{2n-2}^2)\ge
v(\mu_2^{2n-2}\mu_4^{2n-4}\do\mu_{2n-2}^2)$$
hence
$$v(\l_2^{n-1}\l_4^{n-2}\do\l_{2n-2}^1)\ge
v(\mu_2^{n-1}\mu_4^{n-2}\do\mu_{2n-2}^1).$$
For any $z\in E_\Ph$ we denote by $\Th(z)\in\Ph$ the trace of the 
$(4n^2-4n)$-th exterior power of $\ad(z):\End(V_\Ph)@>>>\End(V_\Ph)$. We have
$\Th(\tx)=\prod_{1\le i<j\le n}(\l_{2i}-\l_{2j})^8$. Note that 
$$\align&v(\Th(\tx))=8\sum_{i<j}v(\l_{2i}-\l_{2j})\ge8\sum_{i<j}v(\l_{2i})\\&
=8v(\l_2^{n-1}\l_4^{n-2}\do\l_{2n-2}^1)\ge
v(\mu_2^{n-1}\mu_4^{n-2}\do\mu_{2n-2}^1)\\&=
8\sum_{i<j}v(\mu_{2i})=8\sum_{i<j}v(\mu_{2i}-\mu_{2j})=v\Th(x+\e J).\endalign$$
Moreover if $v\Th(\tx)=v\Th(x+\e J)$ then $v(\l_{2i}-\l_{2j})=v(\l_{2i})$ for
all $i<j$ in $[1,n]$ and $v(\l_1\l_2\do\l_s)=v(\mu_1\mu_2\do\mu_s)$ for 
$s=2,4,\do,2n-2$ hence $v(\l_i^2)=v(\mu_i^2)$ that is, $v(\l_i)=v(\mu_i)$ for 
$i=2,4,\do,2n-2$.

Let $\cs$ be the set of all $\tx$ in $x+\e E_\ph$ such that
$v(\Th(\tx))=v(\Th(x+\e J))$. Then $\cs\ne\em$ (it contains $x+\e J$) and is
"Zariski open" in $x+\e E_\ph$. Let $\tx\in\cs$ and let $\l_1,\do,\l_{2n}$ be
as above. We have 

(b) $v(\l_1)=v(\l_2)=\do=v(\l_{2n_1})=1/n_1$,

$v(\l_{2n_1+1}=v(\l_{2n_1+2})=\do=v(\l_{2n_1+2n_2})=1/n_2$, etc.
\nl
except that when $n_t=1$, $v(\l_{2n}-1)=v(\l_{2n})$ is an integer (or $\iy$)
not necessarily equal to $v(\mu_{2n})$. As in \cite{\KL, p.156} we see that
$$v(\l_i)=1/m\implies\l_i\in\CC[[\e^{1/m}]].$$
(In \cite{\KL, p.156} lines $-7,-8$  one should replace "$t$ not divisible by
$m$" by "$t$ not a divisor of $m$".) This and (b) implies that 
$\tx\in E'_{\Ph,\s}$ where $\s$ contains a product of disjoint cycles of length
$n_1,n_2,\do,n_t$. Thus we have $\cs\sub E'_{\Ph,\s}$. We see that the map 
$\Ps$ in \S14 is bijective in this case.

\subhead 16. Example: $(SO_{2n},SO_{2n-1})$\endsubhead
In this subsection we assume that $\kk=\CC$ and that $(G,K)=(SO(V),SO(U))$ with
$V,U,(,)$ as in \S5 and with $\dim(U)$ odd. Let $\Ph,\ph,\bar\Ph,v$ be as in 
\S15. Let $U_\Ph=\Ph\ot U,U_\ph=\ph\ot U$, $U_{\bar\Ph}=\bar\Ph\ot U$. Let 
$U'_\Ph$ be the set of all $x\in U_\Ph$ such that $(x,x)\ne0$. Let $U'_{\Ph,1}$
(resp. $U'_{\Ph,-1}$) be the set of all $x\in U'_\Ph$ such that 
$v((x,x))\in2\ZZ$ (resp. $v((x,x))\in2\ZZ+1$. We have 
$U'_\Ph=U'_{\Ph,1}\sqc U'_{\Ph,-1}$.

Let $U_0=\{x\in U;(x,x)=0\}$. Let $x\in U_0$. Let 
$X=u_0+\e u_1+\e^2 u_1+\do\in U_\ph$ where $u_0,u_1,\do$ are in $U$. We have 

$(x+\e X,x+\e X)=2(x,u_0)\e+(2(x,u_1)+(u_0,u_0))\e^2+a\e^3$
\nl
with $a\in\ph$. If
$x=0$ let $\cs$ be the set of all $x+\e X$ (with $X$ varying as above) such 
that $(u_0,u_0)\ne0$. 

If $x\in U_0-\{0\}$ let $\cs$ be the set of all $x+\e X$ (with $X$ varying as 
above) such that $(x,u_0)\ne0$. In any case $\cs$ is a nonempty "Zarisky open"
subset of $x+\e U_\ph$. If $x=0$ we have $\cs\sub U'_{\Ph,1}$. If 
$x\in U_0-\{0\}$ we have $\cs\sub U'_{\Ph,-1}$.

We see that the map $\Ps$ in \S14 is bijective in this case.

\subhead 17. Example: $(GL_4,GL_2\T GL_2)$\endsubhead
In this subsection we assume that $(G,K)=(GL(V'\T V''),GL(V')\T GL(V''))$ where
$V',V''$ are $\kk$-vector spaces of dimension $2$ and $\kk=\CC$. We may 
identify $\fp$ with $\Hom(V',V'')\T\Hom(V'',V')$ with the obvious action of 
$K$. Let $\Ph,\ph,\bar\Ph,v$ be as in \S15. Let $V'_\Ph=\Ph\ot V'$, 
$V''_\Ph=\Ph\ot V'$. In our case $\fp_{\Ph,rs}$ consists of all 
$(A,B)\in\Hom(V'_\Ph,V''_\Ph)\T\Hom(V''_\Ph,V'_\Ph)$ such that 
$AB:V''_\Ph@>>>V''_\Ph$ is a regular semisimple automorphism (or equivalently,
$BA:V'_\Ph@>>>V'_\Ph$ is a regular semisimple automorphism). The group $\cw$ in
\S14 is in our case a dihedral group of order $8$; hence $\un{\cw}$ has five 
elements, say $\g_1,\g_2,\g_3,\g_4,\g_5$. We describe the corresponding 
partition $\fp_{\Ph,rs}=\sqc_{j\in[1,5]}\fp_{\Ph,rs,\g_j}$.

$\fp_{\Ph,rs,\g_j}$ consists of all $(A,B)\in\fp_{\Ph,rs}$ such that the
eigenvalues $\l,\mu$ of $AB$ in $\bar\Ph$ satisfy the following condition:

$j=1$: $\l\in\Ph,\mu\in\Ph$, $v(\l)\in2\ZZ$, $v(\mu)\in2\ZZ$.

$j=2$: $\l\in\Ph,\mu\in\Ph$, $v(\l)\in2\ZZ+1$, $v(\mu)\in2\ZZ+1$.

$j=3$: $\l\in\Ph,\mu\in\Ph$, $v(\l)\in2\ZZ$, $v(\mu)\in2\ZZ+1$.

$j=4$: $v(\l)\in(1/2)+\ZZ$, $v(\mu)\in(1/2)+\ZZ$.

$j=5$: $\l\n\Ph,\mu\n\Ph$, $v(\l)\in\ZZ$, $v(\mu)\in\ZZ$.
\nl
Let $e_1,e_2$ be a basis of $V'$; let $e_3,e_4$ be a basis of $V''$. The 
following elements form a set of representatives for the $K$-orbits on $\cn$.

$N_1: e_1\m e_4, e_2\m e_3, e_3\m0, e_4\m e_2$;

$N_2: e_1\m0, e_2\m e_4, e_3\m e_2, e_4\m e_1$;

$N_3: e_1\m e_4, e_2\m0, e_3\m0, e_4\m e_2$;

$N_4: e_1\m0, e_2\m e_4, e_3\m0, e_4\m e_1$;

$N_5: e_1\m e_4, e_2\m e_3, e_3\m0, e_4\m0$;

$N_6: e_1\m0, e_2\m0, e_3\m e_2, e_4\m e_1$;

$N_7: e_1\m e_4, e_2\m0, e_3\m e_2, e_4\m0$;

$N_8: e_1\m e_4, e_2\m0, e_3\m0, e_4\m0$;

$N_9: e_1\m0, e_2\m0, e_3\m e_1, e_4\m0$;

$N_10=0$.
\nl
For $i\in[1,10]$ we consider $(A_i,B_i)=N_i+\e\x$ where $\x\in\fp_\ph$ is given
by

$e_1\m ae_3+be_4, e_2\m ce_3+de_4, e_3\m xe_1+ye_2, e_4\m ze_1+ue_2$
\nl
and $a,b,c,d,x,y,z,u\in\ph$. Define $a_0,b_0,c_0,d_0,x_0,y_0,z_0,u_0\in\CC$ to
be the constant terms of $a,b,c,d,x,y,z,u$. The characteristic polynomial $P_i$
of $A_iB_i$ is given by

$P_1=X^2-O(\e)X+\e x_0+O(\e^2)$;

$P_2=X^2-O(\e)X+\e a_0+O(\e^2)$;

$P_3=X^2-(\e(z_0+d_0)+O(\e^2))X+\e^2c_0x_0+O(\e^3)$;

$P_4=X^2-(\e(b_0+u_0)+O(\e^2))X+\e^2 a_0y_0+O(\e^3)$;

$P_5=X^2-(\e(z_0+y_0)+O(\e^2))X+\e^2(-x_0u_0+y_0z_0)+O(\e^3)$;

$P_6=X^2-(\e(b_0+c_0)+O(\e^2))X+\e^2(-a_0d_0+b_0c_0)+O(\e^3)$;

$P_7=X^2-(\e(z_0+c_0)+O(\e^2))X+\e^2c_0z_0+O(\e^3)$;

$P_8=X^2-(\e z_0+O(\e^2))X+\e^3c_0(x_0u_0-y_0z_0)+O(\e^4)$;

$P_9=X^2-(\e a_0+O(\e^2))X+\e^3u_0(a_0d_0-b_0c_0)+O(\e^4)$;

$P_{10}=X^2-(\e^2(a_0x_0+b_0z_0+c_0y_0+d_0u_0)+O(\e^3))X
+\e^4(a_0d_0-b_0c_0)(x_0u_0-y_0z_0)+O(\e^5)$,
\nl
where $O(\e^m)$ denotes an element of $\e^m\ph$. We see that:

if $x_0\ne0$ then $N_1+\e\x\in\fp_{\Ph,rs,\g_4}$;

if $a_0\ne0$ then $N_2+\e\x\in\fp_{\Ph,rs,\g_4}$;

if $c_0x_0((z_0+d_0)^2-4c_0x_0)\ne0$, then $N_3+\e\x\in\fp_{\Ph,rs,\g_2}$;

if $a_0y_0((b_0+u_0)^2-4a_0y_0)\ne0$, then $N_4+\e\x\in\fp_{\Ph,rs,\g_2}$;

if $(x_0u_0-y_0z_0)((z_0-y_0)^2+4x_0u_0)\ne0$, then 
$N_5+\e\x\in\fp_{\Ph,rs,\g_2}$;

if $(a_0d_0-b_0c_0)((b_0-c_0)^2+4a_0d_0)\ne0$, then 
$N_6+\e\x\in\fp_{\Ph,rs,\g_2}$;

if $c_0z_0(z_0-c_0)\ne0$ then $N_7+\e\x\in\fp_{\Ph,rs,\g_2}$;

if $c_0z_0(x_0u_0-y_0z_0)\ne0$ then $N_8+\e\x\in\fp_{\Ph,rs,\g_3}$;

if $a_0u_0(a_0d_0-b_0c_0)\ne0$ then $N_9+\e\x\in\fp_{\Ph,rs,\g_3}$;

if $(a_0d_0-b_0c_0)(x_0u_0-y_0z_0)\ne0$ and

$(a_0x_0+b_0z_0+c_0y_0+d_0u_0)^2-4(a_0d_0-b_0c_0)(x_0u_0-y_0z_0)\ne0$,
\nl
then $N_{10}+\e\x\in\fp_{\Ph,rs,\g_1}$. 

Thus we have

$\Ps(N_1)=\Ps(N_2)=\g_4$,

$\Ps(N_3)=\Ps(N_4)=\Ps(N_5)=\Ps(N_6)=\Ps(N_7)=\g_2$,

$\Ps(N_8)=\Ps(N_9)=\g_3$,

$\Ps(N_{10})=\g_1$.

\subhead 18. Final comments\endsubhead
We want to define the notion of symmetric space without any assumption on $p$.
Let $G$ be a connected reductive group over $\kk$ and let $K$ be
a closed connected reductive subgroup of $G$. For simplicity we assume that 
$G,K$ contain a common maximal torus $T$. (A similar definition can be given 
without this assumption.) We have inclusions $R_K\sub R_G\sub X(T)$ where
$X(T)$ is the character group of $T$ and $R_K$ (resp. $R_G$) is the set of 
roots of $K$ (resp. $G$). We say that $(G,K)$ is a symmetric space (of equal 
rank) if the inclusions $R_K\sub R_G\sub X(T)$ are the same as the 
corresponding inclusions for a symmetric space of equal rank in characteristic 
$0$. Thus $K$ does not necessarily come from an involution of $G$. For example
if $p=2$ we can take $G$ of type $E_8$ and $K$ a subgroup of type $D_8$ (such a
subgroup exists but it is not the fixed point set of an involution of $G$). It
would be interesting to see how many of the basic properties of symmetric 
spaces in characteristic $0$ extend to this more general case (including the 
case $p=2$).

\enddocument